\documentclass[twoside,11pt,leqno]{article}
\usepackage{amsfonts}

   \def\H{{\cal H}}

 \def\iso{\textrm{iso}}

\newtheorem{df}{Definition}[section]
\newtheorem{thm}[df]{Theorem} 
\newtheorem{cor}[df]{Corollary} 
\newtheorem{rema}[df] {Remark} 
\def\sfstp{{\hskip-1em}{\bf.}{\hskip1em}}

\def\subject#1{\renewcommand{\thefootnote}{}\footnote
{AMS(MOS) subject classification (2010). Primary: {#1}}}

\def\keywords#1{\renewcommand{\thefootnote}{}\footnote
{Keywords: {#1}}}

\def\enddemo{\qed \endtrivlist} \expandafter\let\csname
enddemo*\endcsname=\enddemo

\def\qedsymbol{\ifmmode\bgroup\else$\bgroup\aftergroup$\fi
\vcenter{\hrule\hbox{\vrule
height.5em\kern.5em\vrule}\hrule}\egroup}
\def\qed{\ifmmode\else\unskip\nobreak\fi\quad\qedsymbol}

\title{\bf  On $n$th roots of normal operators}
\author{\normalsize B.P.~Duggal, I.H.~Kim}
\date{}

\begin{document}

\maketitle \thispagestyle{empty} \vskip-16pt

\subject{Primary47A05, 47A55 Secondary47A80, 47A10.} \keywords{Normal operator, $n$-th root, property $(\beta)_{\epsilon}$, decomposable, quasi-nilpotent part, pole, dominant operator, Weyl and Browder theorems  }
\footnote{}
\begin{abstract}
 For $n$-normal operators $A$ \cite{{AP}, {CN}, {CLTU}}, equivalently $n$-th roots $A$ of normal Hilbert space operators, both $A$ and $A^*$ satisfy the Bishop--Eschmeier--Putinar property $(\beta)_{\epsilon}$, $A$ is decomposable and the quasi-nilpotent part $H_0(A-\lambda)$ of $A$  satisfies $H_0(A-\lambda)^{-1}(0)=(A-\lambda)^{-1}(0)$ for every non-zero complex $\lambda$. $A$ satisfies every Weyl and Browder type theorem, and a sufficient condition for $A$ to be normal is that either $A$ is dominant or $A$ is a class ${\mathcal A}(1,1)$ operator.
\end{abstract}


\section {\sfstp Introduction} Let $B(\H)$ denote the algebra of operators, equivalently bounded linear transformations, on a complex infinite dimensional Hilbert space $\H$ into itself. Every {\em normal operator} $A\in B(\H)$, i.e.,  $A\in B(\H)$ such that $[A^*,A]=A^*A-AA^*=0$, has an $n$th root for every positive integer $n>1$. Thus given a normal $A\in B(\H)$, there exists  $B\in B(\H)$ such that $B^n=A$ (and then $\sigma (B^n)=\sigma(B)^n=\sigma(A)$). A straight forward application of the Putnam-Fuglede commutativity theorem (\cite[Page 103]{Hal}) applied to $[B,B^n]=0$ then implies $[B^*,B^n]=0$. (Conversely, $[B^*,B^n]=0$ implies $B^n$ is normal). Operators  $B\in B(\H)$ satisfying $[B^*,B^n]=0$ have been called $n$-{\em normal}, and a study of the spectral structure of $n$-normal operators, with emphasis on the properties which $B$ inherits from its normal avatar $B^n$, has been carried out in (\cite{AP},~\cite{CN},~\cite{CLTU}).

\medskip

Given  $A\in B(\H)$, let $\sigma(A)\subseteq\angle < \frac{2\pi}{n}$ denote that $\sigma(A)$ is contained in an angle $\angle$, with vertex at the origin, of width less than $\frac{2\pi}{n}$. Assuming $\sigma(B)\subseteq \angle <\frac{2\pi}{n}$ for an $n$-normal operator in $B\in B(\H)$, the authors of (\cite{AP},~\cite{CN},~\cite{CLTU}) prove that $B$ inherits a number of properties from $B^n$, amongst them that $B$ satisfies Bishop-Eschmeier-Putinar property $(\beta)_\epsilon$, $B$ is polaroid (hence also isoloid) and $\lim_{m\rightarrow \infty}\langle x_m, y_m \rangle=0$ for sequences $\{x_m\}, \{y_m\} \subset \H$ of unit vectors such that $\lim_{m\rightarrow \infty}\left\| (B-\lambda)x_m\right\|=0=\lim_{m\rightarrow \infty}\left\| (B-\mu)y_m\right\|$ for distinct scalars $\lambda, \mu \in \sigma_a(B)$. (All our notation is explained in the following section.) That $B$ inherits a property from $B^n$ in many a case has little to do with the normality of $B^n$, but is instead a consequence of the fact that $B^n$ has the property. Thus, if the approximate point spectrum $\sigma_a(B^n)=\sigma_a(B)^n$ of $B^n$ is normal (recall: $\lambda\in\sigma_a(B^n)$ is normal if $\lim_{m\rightarrow \infty}\left\| (B^n-\lambda)x_m\right\|=0$ for a sequence $\{x_m\}\subseteq \H$ of unit vectors  implies $\lim_{m\rightarrow \infty}\left\| (B^n-\lambda)^*x_m\right\|=0$; hyponormal operators, indeed dominant operators, satisfy this property), $\sigma(B)\subseteq \angle <\frac{2\pi}{n}$, and $\{x_m\}$,  $\{y_m\}$ are sequences of unit vectors in $\H$ such that $\lim_{m\rightarrow \infty}\left\| (B^n-\lambda^n)x_m\right\|=0=\lim_{m\rightarrow \infty}\left\| (B^n-\mu^n)y_m\right\|$ for some distinct $\lambda, \mu \in \sigma_a(B)$, then
$$
\lim_{m\rightarrow \infty}\lambda^n\langle x_m, y_m \rangle=\lim_{m\rightarrow \infty}\langle B^nx_m, y_m \rangle=\lim_{m\rightarrow \infty}\langle x_m, {B^*}^ny_m \rangle=\mu^n\lim_{m\rightarrow \infty}\langle x_m, y_m \rangle
$$
implies
$$
(\lambda-\mu)\lim_{m\rightarrow \infty}\langle x_m, y_m \rangle=0 \Longleftrightarrow \lim_{m\rightarrow \infty}\langle x_m, y_m \rangle=0
$$
({\em cf.} \cite[Theorem 2.4]{CN}). It is well known that $w$-hyponormal operators satisfy property $(\beta)_\epsilon$ (\cite{BZ}). If $B^n\in (\beta)_\epsilon$ (i.e., $B^n$ satisfies property $(\beta)_\epsilon$) and $\sigma(B)\subseteq \angle <\frac{2\pi}{n}$, then \cite[Theorem 2.9 and Corollary 2.10]{D1} imply that $B+N \in (\beta)_\epsilon$ for every nilpotent operator $N$ which commutes with $B$ ({\em cf.} \cite[Theorem 3.1]{CLTU}). Again, if $B^n$ is polaroid and $\sigma(B)\subseteq \angle <\frac{2\pi}{n}$, then $B$ is polaroid (hence also, isoloid) (\cite[Theorem 4.1]{DDHZ}). Observe that paranormal operators are polaroid. $N$th roots of normal operators have been studied by a large number of authors (see \cite{St},~\cite{RR},~\cite{CF},~\cite{Em},~\cite{Gi}) and there is a rich body of text available in the literature. Our starting point in this note is that  an $n$-normal operator $B$ considered as an $n$th root of a normal operator  has a well defined structure (\cite[Theorem 3.1]{Gi}). The problem then is that of determining the "normal like" properties which $B$ inherits. We prove in the following that the condition $\sigma(B)\subseteq \angle <\frac{2\pi}{n}$ may be dispensed with in many a case (though not always). Just like normal operators, $n$th roots $B$ have SVEP (the single-valued extension property) everywhere, $\sigma(B)=\sigma_a(B)$, $B$ is polaroid (hence also, isoloid). $B\in (\beta)_\epsilon$ (as also does $B^*$) and (the quasinilpotent part) $H_0(B-\lambda)=(B-\lambda)^{-1}(0)$ at every $\lambda \in \sigma_p(B)$ except for $\lambda=0$ when we have $H_0(B)=B^{-n}(0)$. Again, just as for normal operators,  $B$ satisfies various variants of the classical Weyl's theorem $\sigma(B)\setminus\sigma_w(B)=E_0(B)$ (resp., Browder's theorem $ \sigma(B)\setminus\sigma_w(B)=\Pi_0(B)$). It is proved that dominant and class ${\mathcal A}(1,1)$  operators $B$ are normal.

\section {\sfstp Notation and terminology} Given an operator $S\in B(\H)$, {\em the point spectrum, the approximate point spectrum, the surjectivity spectrum} and {\em the spectrum} of $S$ will be denoted by $\sigma_p(S), \sigma_a(S), \sigma_{su}(S)$ and $\sigma(S)$, respectively. The isolated points of a subset $K$ of $\mathbb C$, the set of complex numbers, will be denoted by iso($K$). An operator $X\in B(\H)$ is a {\em quasi-affinity} if it is injective and has a dense range, and operators $S,T \in B(\H)$ are {\em quasi-similar} if there exist quasi-affinities $X,Y\in B(\H)$ such that $SX=XT$ and $YS=TY$.

\medskip

$S\in B(\H)$ has {\em SVEP, the single-valued extension property}, at a point $\lambda_0\in {\mathbb C}$ if for every open disc $\mathfrak{D}$ centered at $\lambda_0$ the only analytic function $f:\mathfrak{D}\rightarrow \H$ satisfying $(S-\lambda)f(\lambda)= 0$ is the function $f\equiv 0$; $S$ has {\em SVEP} if it has SVEP everywhere in $\mathbb C$. (Here and in the sequel, we write $S-\lambda$ for $S-\lambda I$.) Let, for an open subset $\mathcal{U}$ of $\mathbb C$, $\mathcal{E}(\mathcal{U},\H)$ (resp., $\mathcal{O}(\mathcal{U},\H)$) denote the Fr$\rm\acute{e}$chet space of all infinitely differentiable (resp., analytic) $H$-valued functions on $\mathcal U$ endowed with the topology of uniform convergence of all derivatives (resp., topology of uniform convergence) on compact subsets of $\mathcal{U}$. $S\in B(\H)$ satisfies {\em property} $(\beta)_\epsilon$, $S\in (\beta)_\epsilon$, at $\lambda\in\mathbb C$ if there exists a neighborhood $\mathcal{N}$ of $\lambda$ such that for each subset  $\mathcal{U}$ of $\mathcal{N}$ and sequence $\{f_n\}$ of $H$-valued functions in $\mathcal{E}(\mathcal{U},\H)$,
$$
(S-z)f_n(z)\rightarrow 0 \ {\rm in} \ \mathcal{E}(\mathcal{U},\H) \Longrightarrow f_n(z)\rightarrow 0 \ {\rm in} \ \mathcal{E}(\mathcal{U},\H)
$$
(resp., $S$ satisfies {\em property} $(\beta)$, $S\in (\beta)$, at $\lambda \in \mathbb C$ if there exists an $r>0$ such that, for every open subset $\mathcal{U}$ of the open disc $\mathfrak{D}(\lambda;r)$ of radius $r$ centered at $\lambda$ and sequence $\{f_n\}$ of $\H$-valued functions in $\mathcal{O}(\mathcal{U},\H)$,
$$
(S-z)f_n(z)\rightarrow 0 \ {\rm in} \ \mathcal{O}(\mathcal{U},\H) \Longrightarrow f_n(z)\rightarrow 0 \ {\rm in} \ \mathcal{O}(\mathcal{U},\H)).
$$
The following implications are well known (\cite{EP},~\cite{LN}):
$$
S\in (\beta)_\epsilon \Longrightarrow S\in (\beta) \Longrightarrow S \ {\rm has} \ {\rm SVEP};\ S, S^*\in (\beta) \Longrightarrow S \ {\rm decomposable}.
$$

The {\em ascent} asc$(S-\lambda)$ (resp., {\em descent} dsc$(S-\lambda)$) of $S$ at $\lambda\in\mathbb{ C}$ is the least non-negative integer $p$ such that $(S-\lambda)^{-p}(0)=(S-\lambda)^{-(p+1)}(0)$ (resp., $(S-\lambda)^{p}(\H)=(S-\lambda)^{(p+1)}(\H)$). A point $\lambda\in {\rm iso}\sigma(S)$ (resp., $\lambda\in {\rm iso}\sigma_a(S)$) is a {\em pole} (resp., {\em left pole}) of the resolvent of $S$ if $0<{\rm asc}(S-\lambda)={\rm dsc}(S-\lambda)<\infty$ (resp., there exists a positive integer $p$ such that ${\rm asc}(S-\lambda)=p$ and $(S-\lambda)^{p+1}(\H)$ is closed) (\cite{A}). Let
\begin{eqnarray*}
& &\Pi(S)=\{ \lambda\in {\rm iso}\sigma(S): \lambda \ {\rm is \ a\ pole\ (of\ the\ resolvent )\ of} \ S \};\\
& &\Pi^a(S)=\{ \lambda\in {\rm iso}\sigma_a(S): \lambda \ {\rm is \ a\ left \ pole\ (of\ the\ resolvent )\ of} \ S \}.
\end{eqnarray*}
Then $\Pi(S)\subseteq\Pi^a(S)$, and $\Pi^a(S)=\Pi(S)$ if (and only if) $S^*$ has SVEP at points  $\lambda\in\Pi^a(S)$. We say in the following that the operator $S$ is {\em polaroid} if $\{ \lambda\in {\mathbb C}: \lambda \in {\rm iso}\sigma(S)\}\subseteq \Pi(S)$. Polaroid operators are isoloid (where $S$ is {\em isoloid} if $\{ \lambda\in {\mathbb C}: \lambda \in {\rm iso}\sigma(S)\}\subseteq \sigma_p(S)$). Let $\sigma_x=\sigma$ or $\sigma_a$. The sets $E^x(S)=E(S)$ or $E^a(S)$ and $E_0^x(S)=E_0(S)$ or $E_0^a(S)$ are then defined by
\begin{eqnarray*}
& &E^x(S)=\{ \lambda\in {\rm iso}\sigma_x(S):\lambda\in\sigma_p(S)\}, \ {\rm and}\\
& &E_0^x(S)=\{ \lambda\in {\rm iso}\sigma_x(S):\lambda\in\sigma_p(S), \ {\rm dim}(S-\lambda)^{-1}(0)<\infty\}.
\end{eqnarray*}
It is clear that
$$
\Pi^x(S)\subseteq E^x(S) \ {\rm and}\ \Pi_0^x(S)\subseteq E_0^x(S)
$$
(where $\Pi_0^x(S)=\{ \lambda\in \Pi^x(S):{\rm dim}(S-\lambda)^{-p}(0)<\infty\}$).

The {\em quasi-nilpotent part} $H_0(S)$ and the {\em analytic core} $K(S)$ of $S\in B(\H)$ are the sets
\begin{eqnarray*}
H_0(S)&=&\left\{ x\in \H:\lim_{n\rightarrow \infty}\left\| S^nx\right\|^{\frac{1}{n}}=0\right\}, \ {\rm and}\\
K(S)&=&\{x\in\H:{\rm there\ exists\ a\ sequence}\ \{x_n\}\subset\H\ {\rm and} \ \delta>0\ {\rm for}\\
& & {\rm which}\ x=x_0, Sx_{n+1}=x_n\ {\rm and}\ \|x_n\|\leq\delta^n\|x\| \ {\rm for\ all}\ n=1,2,\cdots \}
\end{eqnarray*}
(\cite{A}). If $\lambda\in {\rm iso}\sigma(S)$, then $\H$ has a direct sum decomposition $\H=H_0(S-\lambda)\oplus K(S-\lambda)$, $S-\lambda{\mid }_{H_0(S-\lambda)}$ is quasinilpotent and $S-\lambda{\mid }_{K(S-\lambda)}$ is invertible. A necessary and sufficient condition for a point $\lambda\in\iso\sigma(S)$ to be a pole of $S$ is that there exist a positive integer $p$ such that $H_0(S-\lambda)=(S-\lambda)^{-p}(0)$.

\bigskip

In the following we shall denote the upper semi-Fredholm, the lower semi-Fredholm and the Fredholm spectrum of $S$ by $\sigma_{usf}(S), \sigma_{lsf}(S)$ and $\sigma_f(S)$; $\sigma_{uw}(S), \sigma_{lw}(S)$ and $\sigma_w(S)$ (resp., $\sigma_{ub}(S), \sigma_{lb}(S)$ and $\sigma_b(S)$) shall denote the upper Weyl, the lower Weyl and the Weyl (resp., the upper Browder, the lower Browder and the Browder) spectrum of $S$. Additionally, we shall denote the upper $B$-Weyl, the lower $B$-Weyl and the $B$-Weyl (resp., the upper $B$-Browder, the lower $B$-Browder and the $B$-Browder) spectrum of $S$ by $\sigma_{ubw}(S), \sigma_{lbw}(S)$ and $\sigma_{bw}(S)$ (resp., $\sigma_{ubb}(S), \sigma_{lbb}(S)$ and $\sigma_{bb}(S)$). We refer the interested reader to the monograph (\cite{A}) for definition, and other relevant information, on these distinguished parts of the spectrum; our interest here in these spectra is at best peripheral.

\section {\sfstp Results.}
Throughout the following, $A\in B(\H)$ shall denote an $n$-normal operator. Considered as an $n$th root of the normal operator $A^n$, $A$ has a direct sum representation
$$
A=\bigoplus_{i=0}^\infty A \mid_{\H_i}=\bigoplus_{i=0}^\infty A_i, \ \H=\bigoplus_{i=0}^\infty \H_i,
$$
where $A_0$ is $n$-nilpotent and $A_i$, for all $i=1,2,\cdots$, is similar to a normal operator $N_i \in B(\H_i)$. Equivalently,
$$
A=B_1\oplus B_0,\ B_0=A_0 \ {\rm and} \ B_1=\bigoplus_{i=1}^\infty A_i,
$$
where $B_0^n=0$ and $B_1$ is quasi-similar to a normal operator $N=\bigoplus_{i=1}^\infty N_i\in B\left(\bigoplus_{i=1}^n \H_i\right)$. Quasi-similar operators preserve SVEP; hence, since the direct sum of operators has SVEP at a point if and only if the summands have SVEP at the point, $A$ and $A^*$ have SVEP (everywhere). Consequently (\cite{A}):
$$
\sigma(A)=\sigma(B_1)\cup \{0\}=\sigma(N)\cup \{0\}=\sigma_a(A)=\sigma_{su}(A),
$$
$$
E^a(A)=E(A),\ E_0^a(A)=E_0(A), \ \Pi^a(A)=\Pi(A),\ \Pi_0^a(A)=\Pi_0(A);
$$
furthermore:
$$
\sigma_f(A)=\sigma_{usf}(A)=\sigma_{lsf}(A)=\sigma_{w}(A)=\sigma_{uw}(A)=\sigma_{lw}(A)=\sigma_{b}(A)=\sigma_{ub}(A)=\sigma_{lb}(A),
$$
$$
\sigma_{bf}(A)=\sigma_{bw}(A)=\sigma_{ubw}(A)=\sigma_{lbw}(A)=\sigma_{bb}(A)=\sigma_{ubb}(A)=\sigma_{lbb}(A).
$$
The point spectrum of a normal operator consists of normal eigenvalues (i.e., the corresponding eigenspaces are reducing): This fails for the operator $A$ (\cite[Remark 2.17]{CN}), and a sufficient condition is that $\sigma(A)\subseteq\angle< \frac{2\pi}{n}$ (for then $(A-\lambda)x=0\Longrightarrow (A^n-\lambda^n)x=0\Longrightarrow ({A^*}^n-{\bar\lambda}^n)x=0\Longleftrightarrow ({A^*}-{\bar\lambda})x=0$).

The polaroid property travels from $A^n$ to $A$, no restriction on $\sigma(A)$. (This would then  imply that $E^a(A)=E(A)=\Pi(A)=\Pi^a(A)$ and $E_0^a(A)=E_0(A)=\Pi_0(A)=\Pi_0^a(A)$.) We start by proving that the quasi-similarity of $B_1$ and $N$ transfers to the Riesz projections $P_{B_1}(\lambda)$ and $P_N(\lambda)$ corresponding to points $\lambda\in {\rm iso}\sigma(B_1)={\rm iso}\sigma(N)$. Let $\Gamma$ be a positively oriented path separating $\lambda$ from $\sigma(B_1)$ and let $X,Y$ be quasi-affinities such that $B_1X=XN$ and $YB_1=NY$. Then, for all $\mu\notin \sigma(B_1)$,
\begin{eqnarray*}
& &P_{B_1}(\lambda)=\frac{1}{2\pi i} \int_\Gamma (\mu-B_1)^{-1}d\mu \Longleftrightarrow YP_{B_1}(\lambda)=Y\left\{\frac{1}{2\pi i} \int_\Gamma (\mu-B_1)^{-1}d\mu\right\}\\
& & \Longleftrightarrow YP_{B_1}(\lambda)=\left\{\frac{1}{2\pi i} \int_\Gamma (\mu-N)^{-1}d\mu\right\}Y=P_N(\lambda)Y.
\end{eqnarray*}
A similar argument proves
$$
P_{B_1}(\lambda)X=XP_N(\lambda).
$$

\begin{thm}\label{thm01} $A$ is polaroid.
\end{thm}

\begin{demo} Continuing with the argument above, the normality of $N$ implies that the range $H_0(N-\lambda)$ of $P_N(\lambda)$ coincides with $(N-\lambda)^{-1}(0)$. Hence $(N-\lambda)P_N(\lambda)=0$, and
\begin{eqnarray*}
& & Y(B_1-\lambda)P_{B_1}(\lambda)=(N-\lambda)YP_{B_1}(\lambda)=(N-\lambda)P_{N}(\lambda)Y=0\\
&\Longrightarrow &  (B_1-\lambda)P_{B_1}(\lambda)=0 \Longleftrightarrow H_0(B_1-\lambda)=(B_1-\lambda)^{-1}(0).
\end{eqnarray*}
Since $\lambda\in {\rm iso}\sigma(B_1)$,
\begin{eqnarray*}
& & \bigoplus_{i=1}^\infty \H_i=H_0(B_1-\lambda)\oplus K(B_1-\lambda)=(B_1-\lambda)^{-1}(0)\oplus K(B_1-\lambda)\\
&\Longrightarrow &\bigoplus_{i=1}^\infty \H_i=(B_1-\lambda)^{-1}(0)\oplus (B_1-\lambda)\bigoplus_{i=1}^\infty \H_i,
\end{eqnarray*}
i.e., $\lambda$ is a (simple) pole. The $n$-nilpotent operator $B_0$ being polaroid, the direct sum $B_0\oplus B_1$ is polaroid (since ${\rm asc}(A-\lambda)\leq{\rm asc}(B_0-\lambda)\oplus{\rm asc}(B_1-\lambda)$ and ${\rm dsc}(A-\lambda)\leq{\rm dsc}(B_0-\lambda)\oplus{\rm dsc}(B_1-\lambda)$ for all $\lambda$ (\cite[Exercise 7, Page 293]{TL} )).
\end{demo}

Theorem \ref{thm01} implies:

\begin{cor}\label{cor01} $A$ is isoloid (i.e., points $\lambda\in {\rm iso}\sigma(A)$ are eigenvalues of $A$).
\end{cor}

More is true and, indeed, Theorem \ref{thm01} is a consequence of the following result which shows that $H_0(A-\lambda)=(A-\lambda)^{-1}(0)$ for all non-zero $\lambda\in \sigma(A)$.

\begin{thm}\label{thm02} $H_0(A-\lambda)=(A-\lambda)^{-1}(0)$ for all non-zero $\lambda\in \sigma(A)$ and $H_0(A)=A^{-n}(0)$. In particular, $A$ is polaroid.
\end{thm}

\begin{demo} Following the same notation as above, the normality of $N$ implies $H_0(N-\lambda)=(N-\lambda)^{-1}(0)$ for all $\lambda\in \sigma(N)$ ($=\sigma(B_1)$). Since
$$
NY=YB_1 \Longleftrightarrow (N-\lambda)Y=Y(B_1-\lambda), \ {\rm all}\ \lambda,
$$
it follows that
$$
\left\|(N-\lambda)^nYx\right\|^{\frac{1}{n}}=\left\|Y(B_1-\lambda)^nx\right\|^{\frac{1}{n}}\leq  \left\|Y\right\|^{\frac{1}{n}}\left\|(B_1-\lambda)^nx\right\|^{\frac{1}{n}}\rightarrow 0 \ {\rm as} \ n\rightarrow \infty
$$
for all $x\in H_0(B_1-\lambda)$. Consequently,
$$
Yx\in H_0(N-\lambda)=(N-\lambda)^{-1}(0)\Longrightarrow Y(B_1-\lambda)x=(N-\lambda)Yx=0 \Longleftrightarrow x\in (B_1-\lambda)^{-1}(0),
$$
and hence
$$
H_0(B_1-\lambda)=(B_1-\lambda)^{-1}(0)
$$
for all $\lambda\in \sigma(B_1)$. Evidently,
$$
H_0(A)=H_0(B_1\oplus B_0)=B_1^{-1}(0)\oplus B_0^{-n}(0)\subseteq A^{-n}(0).
$$
Argue now as in the proof of Theorem \ref{thm01} to prove that $A$ is polaroid.
\end{demo}

The Riesz projection $P_A(\lambda)$ corresponding to points $(0\neq)$ $\lambda\in{\rm iso}\sigma(A)$ are, in general, not self-adjoint. Since $\sigma(A)\subseteq \angle < \frac{2\pi}{n}$ ensures $(A-\lambda)^{-1}(0)\subseteq (A^*-{\bar\lambda})^{-1}(0)$ for all $0\neq\lambda\in\sigma_p(A), \ \sigma(A)\subseteq \angle < \frac{2\pi}{n}$ forces $P_A(\lambda)=P_A(\lambda)^*$ for all $\lambda\neq 0$.

\begin{cor}\label{cor02} If $\sigma(A)\subseteq \angle < \frac{2\pi}{n}$, then the Riesz projection corresponding to non-zero $\lambda\in{\rm iso}\sigma(A)$ is self-adjoint.
\end{cor}

\begin{rema}\label{rema01} Theorem \ref{thm01} and \ref{thm02} generalize corresponding results from \cite{AP},~\cite{CN},~\cite{CLTU} by removing the hypothesis that $\sigma(A)\subseteq \angle < \frac{2\pi}{n}$, and, in the case of Theorem \ref{thm02}, the hypothesis on the points $\lambda$ being isolated in $\sigma(A)$. Recall from \cite[Page 336]{A} that an operator $S\in B(\H)$ is said to have property $Q$ if $H_0(S_\lambda)$ is closed for all $\lambda$: {\em Theorem \ref{thm02} says that the $n$th roots $A$ have property $Q$}. Another proof of Theorem \ref{thm02}, hence also of the fact that the operators $A$ satisfy property $Q$, follows from the argument below proving the subscalarity of $A$.
\end{rema}

Property $(\beta)_\epsilon$ (similarly $(\beta)$) does not travel well under quasi-affinities. Thus $CX=XB$ and $B\in (\beta)_\epsilon$ does not imply $C\in (\beta)_\epsilon$ (see \cite[Remark 2.7]{D1} for an example). However, $C\in (\beta)_\epsilon$ implies $B\in (\beta)_\epsilon$ holds, as the following argument proves. If $\{f_n\}$ is a sequence in $\mathcal{E}(\mathcal{U},\H)$ such that
$$
(B-z)f_n(z)\rightarrow 0 \ {\rm in}\ \mathcal{E}(\mathcal{U},\H),
$$
then
$$
X(B-z)f_n(z)=(C-z)Xf_n(z)\rightarrow 0 \ {\rm in}\ \mathcal{E}(\mathcal{U},\H).
$$
Since $C\in (\beta)_\epsilon$ and $X$ is a quasi-affinity,
$$
Xf_n(z)\rightarrow 0 \ {\rm in}\ \mathcal{E}(\mathcal{U},\H) \Longrightarrow f_n(z)\rightarrow 0 \ {\rm in}\ \mathcal{E}(\mathcal{U},\H).
$$
Thus $B\in (\beta)_\epsilon$.

\begin{thm}\label{thm03} $A$ and $A^*$ satisfy property  $(\beta)_\epsilon$.
\end{thm}

\begin{demo} Recall from \cite[Lemma 2.2]{D1} that a direct sum of operators satisfies $(\beta)_\epsilon$ if and only if the individual operators satisfy $(\beta)_\epsilon$. The operator $A$ being the direct sum $B_1\oplus B_0$, where $B_0, B^*_0$ being nilpotent satisfy $(\beta)_\epsilon$, to prove the theorem it will suffice to prove $B_1,  B_1^*\in (\beta)_\epsilon$. But this is immediate from the argument above, since normal operators $N$ satisfy $N,N^*\in (\beta)_\epsilon$ and since there exist quasi-affinities $X$ and $Y$ in $B\left(\bigoplus_{i=1}^\infty\H_i\right)$ such that $N^*X^*=X^*B_1^*$ and $NY=YB_1$.
 \end{demo}

$A\in (\beta)_\epsilon$ implies $A\in (\beta)$, and  $A, A^*\in (\beta)$ implies $A$ is decomposable (\cite{LN}). Hence:

\begin{cor}\label{cor03} $A$ is decomposable.
\end{cor}

We consider next a sufficient condition for the operator $A$ to be normal. However, before that we point out that the operator $A$ satisfies almost all Weyl and Browder type theorems (\cite{A}) satisfied by normal operators.

\medskip
\noindent{\bf Weyl's theorem} An operator $S\in B(\H)$ satisfies
\begin{eqnarray*}
& &{\rm generalized\ Weyl's \ theorem}, \ S\in {\rm gWt}, \ {\rm if} \ \sigma(S)\setminus \sigma_{Bw}(S)=E(S);\\
& &a-{\rm generalized\ Weyl's \ theorem}, \ S\in a-{\rm gWt}, \ {\rm if} \ \sigma_a(S)\setminus \sigma_{uBw}(S)=E^a(S)
\end{eqnarray*}
(see \cite[Definitions 6.59, 6.81]{A}). Let $S\in {\rm Wt}, S\in a-{\rm Wt}, S\in {\rm gBt}, S\in a-{\rm gBt}, S\in {\rm Bt}$ and $S\in a-{\rm Bt}$ denote, respectively, that
\begin{eqnarray*}
& &S \ {\rm satisfies \  Weyl's \ theorem}: \sigma(S)\setminus \sigma_{w}(S)=E_0(S),\\
& &S \ {\rm satisfies \ a- Weyl's \ theorem}: \sigma_a(S)\setminus \sigma_{aw}(S)=E_0^a(S),\\
& &S \ {\rm satisfies \ generalized\ Browder's \ theorem}: \sigma(S)\setminus \sigma_{Bw}(S)=\Pi(S),\\
& &S \ {\rm satisfies \ generalized\ a-Browder's \ theorem}: \sigma_a(S)\setminus \sigma_{uBw}(S)=\Pi^a(S),\\
& &S \ {\rm satisfies \ Browder's \ theorem}: \sigma(S)\setminus \sigma_{w}(S)=\Pi_0(S),\\
& &S \ {\rm satisfies \ a-Browder's \ theorem}: \sigma_a(S)\setminus \sigma_{aw}(S)=\Pi_0^a(S),\\
\end{eqnarray*}
(see \cite[Chapter 6]{A}). The following implications are well known (\cite[Chapters 5, 6]{A}):
$$
S\in a-{\rm gWt}\Longrightarrow \left\{ \begin{array}{l}
S\in a-{\rm Wt}\\
S\in {\rm gWt}
\end{array}\right.
\Longrightarrow S\in {\rm Wt}\Longrightarrow S\in {\rm Bt},
$$
$$
S\in a-{\rm gWt}\Longrightarrow \left\{ \begin{array}{l}
S\in a-{\rm Wt}\\
S\in a-{\rm gBt}
\end{array}\right.
\Longrightarrow S\in a-{\rm Bt}\Longrightarrow S\in {\rm Bt},
$$
$$
S\in a-{\rm gBt}\Longleftrightarrow S\in a-{\rm Bt}, \ S\in {\rm gBt}\Longleftrightarrow S\in {\rm Bt}.
$$
$A$ has SVEP (guarantees $A\in a-{\rm gBt}$ (\cite[Therem 5.37]{A})) and $\sigma(A)=\sigma_a(A)$ guarantee the equivalence of a-gBt and gBt (hence also of a-gBt with a-Bt and Bt) for $A$. The fact that $A$ is polaroid and $\sigma(A)=\sigma_a(A)$ guarantees also that $E(A)=E^a(A)=\Pi^a(A)=\Pi(a)$ (and $E_0(A)=E_0^a(A)=\Pi_0^a(A)=\Pi_0(a)$). Hence all Weyl's theorems (listed above) are equivalent for $A$ and :

\begin{thm}\label{thm04} $A\in a-{\rm gWt}$
\end{thm}

\medskip

\noindent{\bf Normal $A$.}
For the operator $A=B_1\oplus B_0$ to have any chance of being a normal operator, it is necessary that (either $B_0$ is missing, or) $B_0=0$. The hypothesis ($B_0$ is missing, or) $B_0=0$ is, however, in no way sufficient to ensure the normality of $A$. Additional hypotheses are required. An operator $S\in B(\H)$ is said to be {\em dominant} (resp., {\em class} ${\mathcal A}(1,1)$) if to every complex $\lambda$ there corresponds a real number $M_\lambda >0$ such that $\left\|(S-\lambda)^*x\right\| \leq M_\lambda \left\|(S-\lambda)x\right\|$ for all $x\in \H$ (resp., $|S|^2\leq \left|S^2\right|$)  (\cite{SW},~\cite{IY}). Recall from \cite[Lemma 2.1]{DDJ} that if a dominant or class ${\mathcal A}(1,1)$ operator $A\in B(\H)$ is a square root of a normal operator, then $A$ is normal. The following theorem, which uses an argument different from that used in \cite{DDJ}, proves that this result extends to $n$th roots $A$.

\begin{thm}\label{thm05} Dominant or ${\mathcal A}(1,1)$ $n$th roots of a normal operator in $B(\H)$ are normal.
\end{thm}

\begin{demo} Recall that the eigenvalues of a dominant operator are normal (i.e., they are simple and the corresponding eigenspace is reducing). Hence if our $n$th root of $A=B_1\oplus B_0$ is dominant, then $A=B_1\oplus 0$ is a dominant operator which satisfies
$$
A\left(Y\oplus I\mid_{\H_o}\right)=\left(Y\oplus I\mid_{\H_o}\right)(N\oplus 0).
$$
The operator $N\oplus 0$ being normal and the operator $Y\oplus I\mid_{\H_o}$ being a quasi-affinity it follows from \cite{SW},~\cite{D2} that $A$ is normal (and unitarily equivalent to $N\oplus 0$). We consider next $A\in{\mathcal A}(1,1)$.

It is well known that ${\mathcal A}(1,1)$ operators have ascent less than or equal to one. (Indeed, operators $S\in{\mathcal A}(1,1)$ are {\em paranormal}: $\|Sx\|^2 \leq \left\|S^2x\right\|\|x\|$ for all $x\in \H$, hence ${\rm asc}(S)\leq 1.$) Hence if $A=B_1\oplus B_0\in {\mathcal A}(1,1)$, then $B_0=0$ and $A\in B\left(A^{-1}(0)\oplus {A^{-1}}(0)^\perp\right)$ has an upper triangular matrix representation
$$
A=\left(\begin{array}{clcr}0 & A_{12}\\ 0 & A_{22}\end{array}\right).
$$
Let $N_1=N\oplus 0\mid_{\H_0}$ have the represenation
$$
N_1=0\oplus N_{22}\in B\left(N_1^{-1}(0)\oplus {N_1^{-1}(0)}^\perp\right),
$$
and let $Y_1=Y\oplus I\mid_{\H_0}\in B\left(N_1^{-1}(0)\oplus {N_1^{-1}(0)}^\perp, A^{-1}(0)\oplus {A^{-1}(0)}^\perp\right)$ have the corresponding matrix representation
$$
Y_1=\left[Y_{ij}\right]_{i,j=1}^2.
$$
Then, given that $Y$ is a quqsi-affinity satisfying $B_1Y=YN$, $Y_1$ is a quasi-affinity such that $AY_1=Y_1N_1$. Consequently, $A_{22}Y_{21}=0$. The operator $A_{22}$ being injective, we must have $Y_{21}=0$ (and then $Y_{11}$ is injective and $Y_{22}$ has a dense range). The operator $A$ being an $n$th root of a normal operator, $A^n$ is normal. Applying the Putnam-Fuglede commutativity theorem to ($AY_1=Y_1N_1\Longrightarrow$) $A^nY_1=Y_1N_1^n$, it follows that ${A^*}^nY_1=Y_1{N_1^*}^n$, and hence $Y_{12}{N_{22}^*}^n=0$. Since the normal operator ${N_{22}^*}^n$ has a dense range, $Y_{12}=0$ (which than implies that $Y_{11}$ and $Y_{22}$ are quasi-affinities). But then $A_{22}^*Y_{22}=Y_{22}N_{22}^*$ and  $A_{22}Y_{22}=Y_{22}N_{22}$ imply that $A_{22}$ is quasi-affinity. Hence, since ($A^nY_1=Y_1N_1^n$ implies also that) $A_{12}A_{22}^{n-1}Y_{11}=0$, $A_{12}=0$. Thus $A=0\oplus A_{22}$, where $A_{22}\in {\mathcal A}(1,1),\ A_{22}^{-1}(0)=\{0\}$ and $A_{22}Y_{22}=Y_{22}N_{22}$. Applying Proposition 2.5 and Lemma 2.2 of \cite{DDJ}, it follows that $A_{22}$ and $N_{22}$ are (unitarily equivalent) normal operators. Conclusion: $A=0\oplus A_{22}$ is a normal $n$th root.
\end{demo}


\vskip10pt \noindent\normalsize\rm B.P. Duggal, 8 Redwood Grove, London W5 4SZ, England (U.K.).\\
\noindent\normalsize \tt e-mail: bpduggal@yahoo.co.uk

\vskip6pt\noindent \noindent\normalsize\rm I. H. Kim, Department of
Mathematics, Incheon National University, Incheon,  22012, Korea.\\
\noindent\normalsize \tt e-mail: ihkim@inu.ac.kr

\end{document}